\documentclass[10pt,twoside]{article}
\usepackage{amssymb}
\usepackage{Latex-document}

\def\C{{\Bbb C}}
\def\R{{\Bbb R}}
\def\co{{\cal O }}

\def\d#1#2{\frac{\displaystyle #1}{\displaystyle #2}}
\def\ov{\overline}

\def\bc{\begin{center}}
\def\ec{\end{center}}

\def\hang{\hangindent\parindent}
\def\textindent#1{\indent\llap{[#1]\enspace}\ignorespaces}
\def\re{\par\hang\textindent}

\title{\bf Some Results Related to Group Actions\vskip -2mm in Several Complex Variables\vskip 6mm}
\author{Xiangyu Zhou\vspace*{-0.5cm}\thanks{Institute of Mathematics, AMSS, Chinese Academy of Sciences, Beijing;
Department of Mathematics, Zhejiang University, Hangzhou, China.
E-mail: xyzhou@math08.math.ac.cn}}
\date{\vspace{-8mm}}

\begin{document}

\maketitle

\thispagestyle{first} \setcounter{page}{743}

\begin{abstract}

\vskip 3mm

In this talk, we'll present some recent results related to group actions in several complex variables. We'll not
aim at giving a complete survey about the topic but giving some our own results and related ones.

We'll divide the results into two cases: compact and noncompact transformation groups. We emphasize some essential
differences between the two cases. In the compact case, we'll mention some results about schlichtness of envelopes
of holomorphy and compactness of automorphism groups of some invariant domains. In the noncompact case, we'll
present our solution of the longstanding problem -- the so-called extended future tube conjecture which asserts
that the extended future tube is a domain of holomorphy. Invariant version of Cartan's lemma about extension of
holomorphic functions from the  subvarities in the sense of group actions will be also mentioned.

\vskip 4.5mm


\end{abstract}

\vskip 12mm

\section{Fundamentals of several complex variables}

\vskip -5mm \hspace{5mm}

About one century ago, Hartogs discovered that there exist some domains in several complex variables on which any
holomorphic functions can be extended to larger domains, being different with one complex variable. This causes a
basic concept -- domain of holomorphy.

{\bf Definition. } {\it A domain of holomorphy in $\C^ n$ is a
domain on which there exists a holomorphic function which can't be
extended holomorphically across any boundary points.}

A domain in $\C^ n$ is called holomorphically convex, if given any
infinite discrete point sequence $z_k$ there exists a holomorphic
function $f$ s.t. $f(z_k)$ is unbounded (or $|f(x_v)|\to
+\infty$). Consequently, there exists a holomorphic function which
tends to $+\infty$ at the boundary. By Cartan-Thullen's theorem, a
domain in $\C^ n$ is a domain of holomorphy if and only if the
domain is Stein, i.e., holomorphically convex.

{\bf Definition. } {\it A function $\varphi$ with value in
$[-\infty, +\infty)$ on the domain $D$ in $\C^ n$ is called
plurisubharmonic (p.s.h.): if (i) $\varphi$ is upper
semicontinuous (i.e., $\{\varphi<c\}$ is open for each $ c \in
\R$, or equivalently $\ov {\lim}_{z\to z_0}\varphi (z) $ $\leq
\varphi (z_0)$ for $ z_0 \in D$); (ii) for each complex line $L:
= \{z_0+tr : z_0\in D\}, \varphi |_{L\cap D}$ is subharmonic
w.r.t. one complex variable $t$.}

An equivalent definition in the sense of distributions is that $ i
\partial \ov
\partial \varphi$ is a positive current; in particular, when
$\varphi$ is $C^2$, this means Levi form $\left(\d {\partial ^ 2
\varphi} {
\partial z_i \partial \ov z_j}\right)\geq 0$ everywhere. In other
words, $ dJ d \varphi \geq 0$, where $J$ is the complex structure.
(If $ i \partial  \ov \partial \varphi >0,$ then $\varphi$ is
called strictly p.s.h.)

{\bf Example. } For a bounded domain or a domain biholomorphic to
a bounded domain, the Bergman kernel $K (z,\ov z)$ is strictly
p.s.h..

A pseudoconvex domain in $\C^ n$ is a domain on which there exists
a p.s.h. function which tends to $+\infty$ at the boundary. It's
easy to see that a holomorphical convex domain is pseudoconvex,
since $|f|^2$ is plurisubharmonic function where $f$ is given in
the consequence of the definition of a Stein domain.

A deep characterization of a domain of holomorphy is given by a
solution to Levi problem which is the converse of the above
statement.

{\bf Fact.} A domain $D$ in $\C^n$ is a domain of holomorphy if
and only if the domain is pseudoconvex.

A natural corresponding concept of the domain of holomorphy in the
setting of complex manifolds (complex spaces) is the so-called
Stein manifold (Stein space), which is defined as a
holomorphically convex and holomorphically separable complex
manifold (space) . A complex manifold (or space with finite
embedding dimension) is Stein if and only if it is a closed
complex submanifold (or subvariety) in some $\C^n,$ and
 if and only if there exists a strictly
p.s.h. exhaustion function which is $\R$-valued (i.e., the value
$-\infty$ is not allowed). A complex reductive Lie group, in
particular a complex semisimple Lie group, is a Stein manifold.

We know that a domain of holomorphy or a Stein manifold are
defined by special holomorphic functions which are usually hard to
construct in several complex variables. However, a pseudoconvex
domain is defined by a special p.s.h. function which is a real
function and then relatively easy to construct. Construction of
various holomorphic objects in several complex variables and
complex geometry is a fundamental and difficult problem. An
important philosophy here is reducing the construction of
holomorphic functions to the construction of plurisubharmonic
functions, because of the solution of Levi problem and H$\ddot
{\mbox o}$rmander's $L^2$ estimates for $\bar {\partial}$ and
other results.

\section{Group actions in several complex variables}

\vskip -5mm \hspace{5mm}

{\bf Definition. } {\it A group action of the group $G$ on a set
$X$ is given by a mapping $\varphi: G\times X\to X$ satisfying the
following: 1) $e\cdot x=x$, 2) $(ab)\cdot x)=a\cdot (b\cdot x)$,
where $e$ is the identity of the group, $a,b,\in G, x\in X$,
$a\cdot x:=\varphi(a,x)$.}

A group action on a set can be restricted on various cases. When
the set is a topological space and the group is a topological
group, the action is continuous, then one gets a topological
transformation group; when the space is a metric space, the
transformation preserves the metric, then one gets a motion group;
when the set is a differentiable manifold and the group is a Lie
group, the action is differentiable,  then one gets a Lie
transformation group; when the set is a vector space, the
transformation preserves the vector space structure, then one gets
a linear transformation group;  when the set is an algebraic
variety (or a scheme), the group is an algebraic group, and  the
action is algebraic, one gets an algebraic transformation group;
when the set is a complex space, the transformation is
holomorphic, and the action is real analytic, then one gets a
(real) holomorphic transformation group (note that in this case,
if the action is continuous then it is also real analytic); if the
set is a complex space, the group is a complex Lie group, and the
the action is holomorphic, then one gets a complex (holomorphic)
transformation group.

In this talk, we're mainly concerned with the last case. We
consider a complex Lie group $G_{\C}$ with a real form $G_{\R}$
acting holomorphically on a complex manifold (also called
holomorphic $G_{\C}$- manifold) and a $G_{\R}$-invariant domain.
It's known that a complex reductive Lie group has a unique maximal
compact subgroup up to conjugate as its real form, but it also has
many noncompact real forms.

A group action on a set can be regarded as a representation of
the group on the whole group of transformations. An effective
group action means the representation is faithful, so it
corresponds to a (closed) subgroup of the whole transformation
group.

Actually, many domains in several complex variables such as
Hartogs, circular, Reinhardt and tube domains can be formulated
in the setting of group actions.

{\bf Examples. } a) Hartogs and circular domains: consider the
Hartogs action of $\C^{\ast}$ with the real form $S^1$ on $\C^ n$:
$\C^{\ast}\times \C^ n\to \C^n$ given by $(t, (z_1,\cdots,z_n))\to
(tz_1,z_2,\cdots,z_n)$, then Hartogs domain is $S^1$-invariant
domain;  consider the circular action of $\C^{\ast}$ with the real
form $S^1$ on $\C^ n$: $\C^{\ast}\times \C^ n\to \C^n$ given by
$(t, (z_1,\cdots,z_n))\to (tz_1,tz_2,\cdots,tz_n)$, then circular
domain is $S^1$-invariant domain.

b) Reinhardt domains: consider the Reinhardt action of
$(\C^{\ast})^n$ on $\C^ n$ given by $$((t_1,\cdots,t_n),
(z_1,\cdots,z_n))\to (t_1z_1,\cdots,t_nz_n),$$ then Reinhardt
domain is $(S^1)^n$-invariant domain. One can similarly defines
matrix Reinhardt domains

c) tube domains: consider the action of $\R^n$ on $\C^ n$ given by
$(r,z)\to r+z$, then $\R^n$-invariant domain is tube domain.

d) future tube: let $M^ 4$ be the Minkowski space with the Lorentz
metric: $ x\cdot y = x_0 y_0- x_1 y _1- x_2 y _2 - x_3 y _3,$
where  $ x=(x_0, x_1, x_2, x_3), y =(y_0, y _1, y _2, y _3) \in R^
4 ;$ let  $V^+ $ and $V^- = - V^+$ be the future and past light
cones in $R^ 4$ respectively, i.e. $V^\pm =\{y \in M: y^ 2>0, \pm
y _0
>0\}$, the corresponding tube domains $\tau^\pm =T^{V^\pm}= R^ 4
+i  V^\pm$ in  $\C^ 4$ are called future and past tubes; let $L$
be the Lorentz group, i.e. $L=O(1,3),$ $L$ has four connected
components, denote the identity component of $L$ by $L^\uparrow
_+,$ which is called the restricted Lorentz group, i.e.
$L^\uparrow _+ = SO _+(1,3);$ let $L (\C)$ be the complex Lorentz
group, i.e.$L=O(1, 3, \C)\cong O(4, \C), L(\C)$ has two connected
components, denote the identity component of $L(\C)$ by $L
_+(\C),$ called the proper complex Lorentz group which has the
restricted Lorentz group as its real form. Considering the linear
action of $L _+(\C)$ on $\C^ 4$, the future (or past) tube is
$L^\uparrow _+$-invariant.

Denote the $N$-point future tube by $ \tau^\pm  _N =\tau^\pm
\times \cdots \times \tau^\pm$ $N$-times, let $L _+ (\C)$ act
diagonally on $\C^{4N},$ i.e. for $ z=(z^{(1)}, \cdots,
z^{(N)})\in \C^{4N}, \wedge z=(\wedge z ^{(1)}, \cdots, \wedge z
^{(N)})$ where  $ \wedge \in L  _+(\C),$ then $ \tau^\pm  _N$ is
$L^\uparrow _+,$-invariant.

e) matrix Reinhardt domains: let $\C^ n[m\times m]=\{(Z_
1,\cdots,Z_ n)$: $Z_ j\in \C [m\times m]\}$ be the space of
$n$-tuples of $m\times m$ matrices. A domain $D\subset \ C^
n[m\times m]$ is called matrix Reinhardt if it is invariant under
the diagonal ${\rm U}(m)\times {\rm U}(m)$ action $(U,V)¡¤(Z_
1,\cdots,Z_ n)\mapsto (UZ_ 1V,\cdots,UZ_ nV)$. These domains are
the usual Reinhardt domains in the case $m=1$. Diag$(D)$ is
defined as the intersection of $D$ with the diagonal matrices $(Z_
1,\cdots,Z_ n)\in \C ^n[m\times m]$

{\bf Slice theory}

When $G$ is a Lie transformation group properly acting on a smooth
manifold $X$ (e.g. when $G$ is compact), one has a satisfactory
slice theory about the structure of a neighborhood of an orbit.
This theory was extended to the case of an affine reductive group
action regularly on an affine variety by D. Luna ([20]) and the
case of a complex reductive Lie group $G$ action holomorphically
on a Stein space $X$ by Snow ([27]). In these cases, the structure
of a neighborhood of a closed orbit is finely determined. We state
the result for reduced Stein spaces. Let $ G\cdot x$ be a closed
orbit, then there exists a locally closed $G_x$-invariant Stein
subspace $B$ containing $x$ s.t. the natural map from the
homogeneous fiber bundle $G\times _{G_x} B$ over $G/G_x\cong G
\cdot x$ is biholomorphic onto a $\pi$-saturated open Stein subset
of  $X$, where $\pi: X\to X// G$ is the categorical quotient (or
GIT quotient) which exists as a Stein space. Here $B$ is called a
slice at $x$. The slice $B$ is transversal to the closed orbit
$G\cdot x$. When $X$ is regular at $x$, then $B$ can be chosen to
be regular.

As a consequence of the slice theorem, one has a stratification of
the categorical quotient $X//G$ at least when $X$ is a Stein
manifold. The stratum with maximal dimension is Zariski open in
$X//G$ and is contained in the regular part of $X//G.$ This is
called principal stratum. The inverse of the principal stratum
under $\pi : X\to X//G$ consists of all $G$-closed orbits
satisfying that they are of maximal dimension $k$ among the
dimensions of all $G$-closed orbits and their corresponding
isotropy groups are of minimum number of components. Such orbits
are called principal closed orbits, and the corresponding isotropy
groups are called principal. When $k=\dim G$, then $X$ is called
having FPIG.

\section{Some results on compact holomorphic transformation groups}

\vskip -5mm \hspace{5mm}

The relationship between orbit connectedness, orbit convexity, and
holomorphical convexity goes back to the beginning of this
century, when several complex variables was born. Due to Hartogs,
Reinhardt, H.Cartan and others, one already knew some classical
relations between completeness, logarithmic convexity and
holomorphical convexity for circular domains, Hartogs domains,
and Reinhardt domains. The orbit connectedness and orbit
convexity are defined in a general setting (for arbitrary compact
connected Lie group), which correspond to completeness and
logarithmic convexity when one restricts to the above domains.

There are some fundamental relationships between orbit
connectedness and orbit convexity with holomorphically convexity
and envelope of holomorphy for invariant domains.

{\bf Definition. } {\it Let $G_{\C}$ be a connected complex Lie
group, $G_{\R}$ be a connected closed real form of $G_{\C}.$ Let
$X$ be a holomorphic $G_{\C}$-space, $D\subset X$ be a
$G_{\R}$-invariant set, we call $D$ orbit connected, if for $
b_z:G_{\C}\to X, g \mapsto g \cdot z, b^{-1}_z (D)$ is connected
for each $ z\in  D.$ When $(G_{\C},G_{\R})$ is a geodesic convex
pair(i.e. the map $\rm{Lie}(G_{\R})\times G_{\R}
\ni(v,g)\rightarrow\exp(iv)¡¤g\in G_{\C}$ is a homeomorphism, cf.
[3]), $D$ is called orbit convex if for each $z\in D,$ and $v\in
i\rm{Lie}(G_{\R})$ s.t. $\exp(v)\in b^{-1}_z (D)$ it follows
$\exp(tv)\in b^{-1}_z (D)$ for all $t\in [0,1]$.}

Roughly speaking, orbit connectedness means that $G_{\C}x\cap D$
is connected for every $x\in D$.

One has known for a long time that the envelope of holomorphy of a
domain in $\C ^n$ (or more general a Riemann domain over $\C ^n$)
exists uniquely as a Riemann domain over $\C ^n$. There is a
difficult problem of univalence: When is the envelope of
holomorphy of a domain in $\C ^ n$ itself a domain in $\C ^n$? We
have the following criteria for the univalence of the envelope of
holomorphy for certain invariant domains:

{\bf Theorem 1 ([36]).} {\it Let $X$ be a Stein manifold, $K^{\C}$
be a complex reductive Lie group holomorphically acting on  $X$,
where $K$ is a connected compact Lie group and $K^{\C}$ be its
universal complexification. Let $D\subset X$ be a $K$-invariant
orbit connected domain. Then the envelope of holomorphy $E(D)$ of
$D$ is schlicht and orbit convex if and only if the envelope of
holomorphy $E(K^{\C}\cdot D)$ of $K^{\C}\cdot D$ is schlicht.
Furthermore, in this case, $E(K^{\C}\cdot D)=K^{\C}\cdot E(D).$}

When $K=S^1$ and the action is circular (or $\alpha$-circular) and
Hartogs, the corresponding concepts of orbit connectedness for
such domains were introduced separately and the above results were
obtained and stated separately by Casadio Tarabushi and Trapani in
[1,2].

When $K=(S^1)^n$ and the action is Reinhardt, the result is well
known as a classical result about Reinhardt domain which asserts
that any Reinhardt domain in $(\C^{\ast})^n$ has schlicht
envelope of holomorphy.

Some other results were also included in the above theorem. So our
result can also be regarded as an extension of their results and
the classical result on Reinhardt domains in a unified way.

In the proof, a theorem due to Harish-Chandra on the infinite
dimensional representation of Lie groups plays an important role.

We also give some examples of orbit connected domains. Let
$X=K^{\C}/L^{\C},$ the action of $K^{\C}$ on $X$ be given by the
left translations. When $L$ is connected or $(K,L)$ is a
symmetric pair, then any $K$-invariant domain is orbit connected.
The simplest example is Reinhardt domains in $(\C^\ast)^n.$

The origin of orbit connectedness could at least go back to [28].

{\bf Example. } A theorem of V.Bargmann, D. Hall and A.S. Wightman (cf. Wightman [32], Jost [12],
Streater-Wightman [28]) asserts that $\tau^+ _N$ is orbit connected.

We also consider the homogeneous embeddings of $K^{\C}/L^{\C}$.
Let $X$ be a smooth homogeneous space embedding of
$K^{\C}/L^{\C}$, $D\subset X$ be a $K$-domain. Assume that $L$ is
connected or $(K,L)$ is a symmetric pair. Then $E(D)$ is schlicht
and orbit convex. In particular, every matrix Reinhardt domain of
holomorphy $D$ is orbit convex.  Since an orbit convex matrix
Reinhardt domain has a path connected Diag$(D)$, so a matrix
Reinhardt domain of holomorphy has a connected Diag$(D)$.

{\bf Theorem 2([37]). } {\it Let $K$ be a connected compact Lie
group, $L$ be a closed (not necessarily connected) subgroup of
$K$. Let $K^{\C}$ and $L^{\C}$ be respectively universal
complexification of $K$ and $L$. Suppose that $D$ is
$K$-invariant relatively compact domain in $K^{\C}/L^{\C}$ (Here
the action of $K^{\C}$ is given by left translations). Then (i)
$\hbox{Aut}(D)$ is a compact Lie group; (ii)  Any proper
holomorphic self-mapping of $D$ is biholomorphic if $K$ is
semisimple or a direct product of a semisimple compact Lie group
and a compact torus.}

By a result of Matsushima, $K^{\C}/L^{\C}$ is a Stein manifold
which is a holomorphic  $K^{\C}$ - manifold w.r.t. left
translation action.

The motivations of the present work are two-folds: the result (i)
is to extend a main result of [4], where  the same result was
obtained by requiring a restrictive condition that $(K,L)$ is a
symmetric pair,i.e., $K/L$ is a compact Riemannian symmetric
space; the result (ii) is to extend a classical result which
asserts that proper self mapping of the relatively compact
Reinhardt domains in $(\C^\ast)^n$ is biholomorphic.

The proof is involved with many famous results such as Mostow
decomposition theorem, H. Cartan's theorem about compactness of
automorphism groups, Andreotti-Frankel's theorem on homology
group of a Stein manifold, the holomorphic version of de Rham's
theorem on a Stein manifold, a result of Milnor's about CW
complex, a result from iteration theory, Poincar\'e duality
theorem, degree theory for proper mappings, covering lifting
existence theorem,  and a result about compact semisimple Lie
groups et al.

\section{Extended future tube conjecture}

\vskip -5mm \hspace{5mm}

Let's keep the notation in Example d of the section 2. The set
$\tau '_N: = \{\wedge z: z\in \tau^+ _N, \wedge \in L_+(\C)\}$ is
called the extended future tube.

The extended future tube conjecture, which arose naturally from
axiomatic quantum field theory at the end of 1950's, asserts that
the extended future tube $\tau' _N $ is a domain of holomorphy for
$ N\geq 3.$ This conjecture turns out to be very beautiful and
natural. In their papers, Vladimirov and Sergeev said that the
importance of the conjecture is also due to the fact that there
are some assertions in QFT, such as the finite covariance theorem
of Bogoliubov-Vladimirov, proved only assuming that this
conjecture is true.

According to the axiomatic quantum field theory (cf. [12,13,28]),
one may describe physical properties of a quantum system using the
Wightman functions which correspond to holomorphic functions in
$\tau^+_N$ invariant w.r.t. the diagonal action of $L_+^\uparrow$.
This sort of functions have the following extension property.

BHW Theorem (due to Bargman, Hall, and Wightman 1957). An
$L_+^\uparrow$-invariant holomorphic function on $\tau^+_N$ can be
extended to an $L_+(\C)$-invariant holomorphic function on
$\tau'_N$ (cf. [12,13,28]).

In the proof, the orbit connectedness of $\tau^+_N$ play a key
role. With this and Identity Theorem, one can easily define the
invariant holomorphic extension.

So, a natural question arises, i.e., can these holomorphic
functions be extended further? Or, is $\tau'_N$ holomorphic
convex w.r.t. $L_+(\C)$-invariant holomorphic function? After
some argument, this is equivalent to ask if $\tau'_N$ is a domain
of holomorphy.

Streater's theorem. A holomorphic function on the Dyson domain
$\tau^+_N\cup \tau^-_N \cup J$ (where $J: =\tau'_N\cap M^{4N}$ is
the set of Jost points which was proved to exist and characterized
by R. Jost) can be extended to a holomorphic function on $\tau'_N$
(cf. [12,28]).

So, a natural question is to construct the envelope of holomorphy
of the Dyson domain $\tau^+_N \cup \tau^-_N \cup J$ (This
question is mentioned in the article ``Quantum field theory" of
the Russian's great dictionary ``Encyclopedia of Mathematics").
That the extended future tube conjecture holds is equivalent to
that this envelope of holomorphy is exactly the extended future
tube $\tau'_N$.

The conjecture have been mentioned as an open problem in many
classical ([12,28]) and recent references ([11,21-24,28-31]) and
references therein. In [38,39], we found a route to solve the
conjecture via Kiselman-Loeb's minimum principle and Luna's slice
theory. Let's recall the minimum principle.

{\bf Minimum principle}

Let $X$ be a complex manifold, $G_{\C}$ a connected complex Lie
group, $G_{\R}$ a connected closed real form of $G_{\C}$. Denote
$\psi : G_{\C} \to G_{\C}/ G_{\R}$, and $p : X \times G_{\C}\to X$
the natural projections.

$G_{\C}$ acts on $X\times G_{\C}$ on the right by:

$$
\begin{array}{rcl}
(X\times G_{\C}) \times G_{\C} &\longrightarrow & X\times G_{\C}\\
((x, g), h) &\longmapsto & (x, g h)
\end{array}
$$

Let $\Omega \subset X \times G_{\C}$ be a right $G_{\R}$-invariant
domain and have connected fibres of $p$; and $u \in C^\infty
(\Omega)$ be a right $G_{\R}$-invariant function. $u$ naturally
induces a smooth function $\dot u(x,\psi(g))$ on $\dot\Omega :
=(id_X, \psi) (\Omega)$.

Suppose that (1) $u$ is p.s.h on $\Omega$, (2) $\forall x \in p
(\Omega), u(x, \cdot)$ is strictly p.s.h. on $\Omega_x = \Omega
\cap p^{-1} (x)$, and (3) $\dot u (x, \cdot)$ is exhaustive on
$\dot\Omega_x =\psi (\Omega_x)$, then the minimum principle
asserts that $v(x) =\inf\limits_{g \in \Omega_x} u (x, g) $  is
$C^\infty$ and p.s.h. on $p(\Omega)$.

{\bf Remark. } C.O. Kiselman in [14] first obtained the minimum
principle when $X=\C^n, G_{\C} =\C^m, G_{\R}=Im \C^m$ , J.J. Loeb
in [18] generalized Kiselman's result to the present general case.

It's easy to construct invariant p.s.h. functions w.r.t. compact
Lie group via ``averaging technique". However, such a technique
doesn't hold again for non compact Lie group.

{\bf Observation. } Let $G$ be a real Lie group which acts on $\C^
n$ linearly. Let $D$ be a Bergman hyperbolic domain which is
$G$-invariant. Then the Bergman kernel $K_D(z,\ov w)$ satisfies
$K _D(z,\ov z)= K _D (g \cdot z, \ov {g \cdot z})$ for $ g\in G,$
when $G$ is compact; when $G$ is semisimple, we have $ K _D(z,
\ov w)= K _D(g\cdot z,\ov {g\cdot w})$.

Brief proof is as follows. Since $G$ linearly act on $\C^n$, one
has a representation $G \to GL(n,\C)$; if $G$ is semisimple, then
the image of $G$ must be in $SL(n,\C);$ if $G$ is compact, the
image of $G$ is in $U(n)$. Using the transformation formula for
the Bergman kernels and noting that the determinant of the
Jacobian of the map $z \to g\cdot z$ is 1 for semisimple case, and
is in $S^1$ for compact case, then we can get the result.

We consider the following question: Let $X$ be a Stein  manifold,
$G_{\C}$ be a connected complex reductive Lie group acting on $X$
s.t. the action is holomorphic, $G_{\R}$ a connected real form of
$G_{\C}$. Let $D\subset X$ be a $G_{\R}$-invariant orbit
connected Stein domain, is $G_{\C}\cdot D$ also Stein?

When $G_{\R}$ is compact, the answer is positive (cf. [22]). This
is a special case of Theorem 1 in the section 3.

The extended future tube conjecture is a special case of the
question, where $X=\C^ {4N}, G_{\C}=L _+(\C), G_{\R}=L^\uparrow
_+, D=\tau^+ _N, G_{\C}\cdot D=\tau '_N$

Consider $X\times G_{\C}\buildrel \rho \over \longrightarrow X,
\rho (x,g)=g^{-1}\cdot x$. Suppose that there is a suitable
$G_{\R}$-invariant s.p.s.h. function $\varphi$ on $D$. We have a
p.s.h.  function $u(x,g)=\varphi (g^{-1}\cdot x)$ on $\Omega =
\rho^{-1}(D).$ Define $\psi (x)= \displaystyle\inf_{g\in
\Omega_x}u(x,g)$ for $ x\in   p (\Omega),$ where $ p: X\times
G_{\C}\to X$ is given by $p(x,g)=x,$ and $\Omega_ x:=\{g\in
G_{\C}: (X,g)\in \Omega\}$.

In order to prove $\psi (x)$ is p.s.h. on $p (\Omega) =
G_{\C}\cdot D,$ we can use the minimum principle due to
Kiselman-Loeb.

{\bf Observation. } $\Omega_x$ is connected if and only if $D$ is
orbit connected.

In order to use the minimum principle, we still need to check two
assumptions: (i) $u(x,\cdot)$ is s.p.s.h. on $\Omega_x;$ (ii)
$\dot u (x,\cdot)$ is exhaustion on $\dot \Omega_x,$ where $\dot u
(x, \psi (g))$ is defined on $\dot \Omega = (id, \psi
)(\Omega)\subset X \times G_{\C}/G_{\R}$ and is induced by
$u,\psi: G_{\C}\to G_{\C}/G_{\R}, \dot \Omega_x = \psi (\Omega_x).
$ Usually speaking, assumption (i) fails on the whole $\Omega$.
However, when  $ X$ has FPIG, then the assumption (i) is fulfilled
on a Zariski open subset of $\Omega$. Let $X':=\{x\in X:  G_{\C} x
\ \hbox{is closed},  \ (G_{\C})_x$ is principal and finite $\},$
then, by the slice theory, $A= X\backslash A'$ is a
$G_{\C}$-invariant analytic subset of $X$. Let $D'=D\cap X',
\Omega':= \rho ^{-1}(D'),$ then the assumption (i) is satisfied on
$\Omega'$. If the assumption (ii) is also satisfied on $\Omega'$,
then we can use the minimum principle on $\Omega'$ and get that
$\psi (x)$ is p.s.h. on  $p (\Omega') = G_{\C} \cdot D \backslash
A$ since $\psi (x)$ is upper semicontinuous on $G_{\C}\cdot D,$ by
the extension theorem for p.s.h. functions, $\psi (x)$ can be
extended to a p.s.h. function on  $G_{\C}\cdot D.$

If we can prove that the extended p.s.h. function is weak
exhaustion, then $G_{\C}\cdot D$ is Stein.

As a consequence of our observations, we deduce that the general
question is true for pseudoconvex pair ($G_{\C},G_{\R}$) (i.e.,
there exists a $G_{\R}$-invariant p.s.h. function on $G_{\C}$
which is exhaustion on $G_{\C}/G_{\R}$(cf.[17]), which include the
case when $G_{\R}$ is compact and nilpotent(cf.[17]). However it's
pity that $(L _+(\C), L^\uparrow _+$) is not a pseudoconvex pair.

In the case of the extended future tube conjecture, we proved that
the assumption (ii) in the minimum principle is satisfied and the
constructed function is weak exhaustion. These are the main
technical difficulties. We overcome them and finished our proof
via a consideration of the matrix form of the conjecture and
explicit calculations based on Hua's work and matrix techniques
([9,19]).

{\bf Theorem [38,39]. } {\it The extended future tube conjecture
is true.}

A.G. Sergeev posed an interesting idea to attack the mentioned
question. He assumed an invariant version of Cartan's lemma: if
$A\subset D $ is a $G_{\R}$-invariant analytic subset, $ f\in
\co(A)^{G_{\R}}$, then there exists an  $F\in \co (D)^{G_{\R}}$
s.t. $F|A=f.$ If this is the case, we can prove that  $\pi (D)$ is
Stein in  $X//G_{\C}$. In order to prove it, it's sufficient to
prove $\pi (D)$ is holomorphically convex. Let $\{y_n\}\subset \pi
(D)$ be an arbitrary discrete set. Then  $\{\pi^{-1} (y_n)\}\cap
D$ is a $G_{\R}$-invariant analytic subset in  $D$. By the
assumption, then there exists a $G_{\R}$-invariant holomorhic
function $F$ on  $D$ s.t. $F|_{\pi^{-1}}(y_n)=n.$ Since $\co(\pi
(D))\cong \co(D)^{G_{\R}},$ then we get a holomorphic function on
$\pi (D)$ which is unbounded on $\{y_n\}.$ This means that $\pi
(D)$ is holomorphically convex, and then  $\pi^{-1}(\pi(D))$ is
also Stein. When $\pi^{-1}(\pi (D))= G_{\C}\cdot D,$ i.e.,
$G_{\C}\cdot  D$ is $\pi$-saturated, then $G_{\C}\cdot D$ is
Stein.

It seems to be hard to prove directly the invariant version of
Cartan's lemma for a noncompact Lie group  $G_{\R}$, although it's
trivially  the case for a compact Lie group. Actually, we have the
following:

{\bf Proposition ([41]). } {\it Suppose, furthermore, $G_{\C}\cdot
D$ is $\pi$-saturated. Then the invariant version of Cartan's
lemma holds if and only if $G_{\C}\cdot D$ is Stein.}

However, we recently observed that it should be possible to
directly give an answer to the above question based on
$L^2$-methods and group actions.

\label{lastpage}

\end{document}